\documentclass[11pt]{amsart}

\usepackage{hyperref}
\usepackage{amsfonts}
\usepackage{amsmath}
\usepackage{amsthm}
\usepackage{latexsym}
\usepackage{amscd}
\usepackage{xypic}
\usepackage[all,cmtip]{xy}
\usepackage{mathrsfs}
\usepackage{graphicx}
\usepackage{enumitem} 
\usepackage{braket}
\usepackage[margin=1.4in]{geometry}
\usepackage{tikz}
\usetikzlibrary{calc,fadings,decorations.pathreplacing}
\usepackage{verbatim}

\allowdisplaybreaks

\numberwithin{equation}{section}

\newtheorem{prop}{Proposition}[section]
\newtheorem{theo}[prop]{Theorem}
\newtheorem{lemm}[prop]{Lemma}
\newtheorem{coro}[prop]{Corollary}

\newtheorem{defi}[prop]{Definition}

\newtheorem{defiprop}[prop]{Definition-Proposition}

\newcommand{\DD}{\mathbb{D}}
\newcommand{\EE}{\mathbb{E}}

\newcommand{\PP}{\mathbb{P}}
\newcommand{\QQ}{\mathbb{Q}}
\newcommand{\RR}{\mathbb{R}}

\newcommand{\cG}{\mathcal G}

\newcommand{\cP}{\mathcal P}

\newcommand\addlabel[1]{\refstepcounter{equation}\tag{\theequation}\label{#1}}

\newcommand{\bangle}[1]{\left\langle #1 \right\rangle}
 
\DeclareMathOperator{\diam}{diam}
\DeclareMathOperator{\Ric}{Ric}

\DeclareMathOperator{\Ent}{Ent}

\DeclareMathOperator{\vol}{vol}

\DeclareMathOperator{\proj}{proj}

\newcommand{\Leb}{\text{\textnormal{Leb}}}
\newcommand{\Id}{\text{Id}}

\DeclareMathOperator{\Var}{Var}

\title[Lack Of Ricci Bounds For The Entropic Measure]{A Lack of Ricci Bounds for the Entropic Measure on Wasserstein Space over the Interval}
\author{Otis Chodosh}
\address{Department of Mathematics, Stanford University, CA 94305-2125, USA}
\email{ochodosh@math.stanford.edu}
\date{\today} 

\begin{document} 
\begin{abstract}
We examine the entropic measure, recently constructed by von Renesse and Sturm, a measure over the metric space of probability measures on the unit interval equipped with the $2$-Wasserstein distance. We show that equipped with this measure, Wasserstein space over the interval does not admit generalized Ricci lower bounds in the entropic displacement convexity sense of Lott-Villani-Sturm. We discuss why this is contrary to what one might expect from heuristic considerations.
\end{abstract}

\maketitle
\tableofcontents

\section{Introduction}

For a compact Riemannian manifold $M$, the space of Borel probability measures on $M$, $\cP(M)$, has very interesting geometric properties when equipped with what is known as the $2$-Wasserstein metric, arising via optimal transport theory. In particular, as discovered by Otto in \cite{Otto:DissapEqn}, it can be formally regarded as an infinite dimensional Riemannian manifold with well understood notions of tangent spaces, geodesics, gradients and more (in fact, many of these these formal notions have very recently been made rigorous in various settings, in addition to Otto's original paper \cite{Otto:DissapEqn}, also see \cite{OttoVillani:GenTalagrand,Lott:WassGeoCalc,Gigli:secondOrder} for discussions of both the formal structure and rigorous results along the same lines). However, there is a notable lack of the notion of a volume measure in this structure. In fact, in \cite[Remark 5.6]{Gigli:secondOrder}, Gigli gives an argument that there should not be a natural volume measure on the space of probability measures, which we discuss in \ref{subsect:heuristics}. If we cannot hope for a natural choice of volume measure ``agreeing with the (formal) Riemannian structure'', we can at least search for measures on $\cP(M)$ with ``nice'' properties. From a geometric point of view, one would be interested in a measure on $\cP(M)$ which would have interesting geometric properties and in particular, in this paper, we will discuss the possibility of lower Ricci curvature bounds on $\cP(M)$. The space $\cP(M)$ is only formally a Riemannian manifold, and it is therefore not a priori clear if this is even a reasonable concept. However, thanks to the work of Lott-Villani-Sturm in \cite{LottVillani:RicciViaTspt,Sturm:MetMeas1,Sturm:MetMeas2}, there is a well defined notion of a metric measure space having lower Ricci bounds: convexity of the Boltzmann entropy functional along Wasserstein geodesics in the space of probability measures on the space of interest. We will give definitions below, but what will be most important about this concept in this paper is that it requires a fixed background measure to make sense of it. Because there is no ``Riemannian volume measure'' on $\cP(M)$, we should instead ask if there are choices which have interesting geometric properties, and in this paper we examine the possibility of one such measure giving the space generalized lower Ricci bounds (which we will see it does not). 

In \cite{RenesseSturm:EntMeasDiff}, von Renesse and Sturm have defined what they call an ``entropic measure'' $\PP^{\beta}_{0}$ on $\cP_{0}:= \cP([0,1])$, probability measures on the unit interval.  In this paper we will show that although it heuristically could be expected to admit lower Ricci bounds, the entropic measure $\PP^{\beta}_{0}$ does not admit \emph{any} lower Ricci bound. We will additionally discuss the possibility of extending our proof to the higher dimensional case for the measure constructed by Sturm in \cite{Sturm:MultEnt} as a generalization of $\PP^{\beta}_{0}$ to arbitrary compact Riemannian manifolds $M$. 

Our main theorem will be (we will provide the necessary definitions in the remainder of the introduction)
\begin{theo}\label{theo:QQbeta-no-ric-geqK}
There is no $K \in \RR, \beta>0$ such that $(\cG_{0},d^{L^{2}},\QQ^{\beta}_{0})$ has generalized $\Ric \geq K$. 
\end{theo}
Because we will see that this space is metric-measure isomorphic to $(\cP_{0},d^{W},\PP^{\beta}_{0})$ (where the metric $d^{W}$ is the 2-Wasserstein distance, as defined below) we will have as an immediate corollary that:
\begin{coro}
There is no $K \in \RR, \beta>0$ such that $(\cP_{0},d^{W},\PP^{\beta}_{0})$ has generalized $\Ric \geq K$.
\end{coro}

\subsection{Optimal Transport and Wasserstein Space} Below, we briefly review the basic concepts of optimal transport, Wasserstein distance, and generalized Ricci bounds through displacement convexity of the entropy functional. For a proper introduction to this topic we highly recommend that the reader refer to Villani's book \cite{Villani:TopOptTspt03}, as well as his monograph \cite{Villani:OptTspt}. In addition to these excellent sources, we remark that the material contained in this paper is a condensed form of the author's essay, \cite{Chodosh:essay}, which contains more precise exposition and references for the introductory material below. 

Let $(X,d)$ be a compact metric space. We will denote the space of probability measures on $X$ by $\cP(X)$. For $k=1,2$, we define $\proj_{k}:X\times X \to X$ to be projection onto the $k$-th factor. For two probability measures $\mu,\nu \in \cP(X)$, we define the set of \emph{admissible transport plans} to be 
\begin{equation} \Pi(\mu,\nu) =\{ \pi \in \cP(X\times X) : (\proj_1)_*\pi = \mu,(\proj_2)_*\pi = \nu\}\end{equation}
where for $k=1,2$, $(\proj_{k})_{*} \pi \in \cP(X)$ is the $k$-th marginal, \[(\proj_{k})_{*}(A) = \pi[(\proj_{k}^{-1}(A)].\] This is certainly nonempty, because the measure $\mu\otimes \nu$ which is defined by \[ \int_{X\times X} f(x,y) d\mu\otimes \nu(x,y) : = \int _{X}\int_{ X} f(x,y) d\mu(x)  d\nu(y)\] for continuous $f\in C(X\times X)$, is clearly in $\Pi(\mu,\nu)$. Loosely speaking, a transport plan is a proposal for how to move the mass of $\mu$ around so as to assemble the distribution of mass prescribed by $\nu$. Given the above definition, we define the $2$-\emph{Wasserstein distance} between $\mu$ and $\nu$ to be 
\begin{equation}\label{eq:wass-dist} d^{W}(\mu,\nu)^{2} := \inf_{\pi\in\Pi(\mu,\nu)}\int_{X\times X} d(x,y)^{2} d\pi(x,y). \end{equation}
One can show that this defines a metric on $\cP(X)$, metrizing the weak-* topology. Furthermore, if $(X,d)$ is a \emph{geodesic space}, meaning that for any two points $x_{0},x_{1}$, there is a continuous curve $\gamma:[0,1] \to X$ with $\gamma(0) = x_{0},\gamma(1) = x_{1}$ and 
\begin{equation} \label{eq:Ldefi}  d(x_{0},x_{1})  = \sup_{a =t_0<t_1<\dots<t_{N+1} = b} \sum_{i=0}^N d(\gamma(t_i),\gamma(t_{i+1})),  \end{equation} 
then it is possible to show that $(\cP(X), d^{W})$ is also a geodesic space.

We now consider an arbitrary  \emph{measured geodesic space} $(X, d, m)$, i.e.\ $(X,d)$ is a metric space which is a geodesic space as defined above, and $m\in \cP(X)$ is a fixed probability measure on $X$. We define the entropy functional $ \Ent(\cdot| m) : \cP(X) \to \overline \RR = \RR \cup \{+ \infty\} $ by  \[ \mu \mapsto \begin{cases}
\int_{X} \rho \log \rho\  dm & \text{ for } \mu \ll m \text {, i.e. } \mu = \rho m\\
+ \infty & \text{ otherwise.}
\end{cases}\]
Intuitively, $\Ent(\mu|m)$ measures the nonuniformity of $\mu$ with respect to $m$.  In some sense, we can think of it as a kind of distance (without symmetry or a triangle inequality), which is ``extensive'' in the sense that if $\mu,\mu',m\in \cP(X)$ then $\Ent(\mu\otimes\mu'|m\otimes m) = \Ent(\mu|m) + \Ent(\mu'|m)$.

\begin{defi}\label{defi:k-convex-ent}
For a measured geodesic space $(X,d,m)$, we say that $\Ent(\cdot | m)$ is weakly a.c.\ $K$-displacement convex on $(X,d,m)$ if for any probability measures $\mu_{0},\mu_{1}\ll m$, there exists a geodesic in $\cP(X)$, $\mu_{t}$ from $\mu_{0}$ to $\mu_{1}$ so that $\Ent(\mu_{t}|m)$ is $K$-convex, in the sense that \begin{equation}\label{eq:k-disp-covex}
\Ent(\mu_{t}|m) \leq t \Ent(\mu_{1}|m)  + (1-t) \Ent(\mu_{0}|m) - \frac K2 t (1-t) d^{W}(\mu_{0},\mu_{1})^{2}. 
\end{equation}
\end{defi}
This definition is motivated by the amazing connection between displacement convexity of the entropy functional with lower Ricci bounds, when $X$ is a Riemannian manifold. 
The following theorem was first proven in this generality by von Renesse-Sturm in \cite[Theorem 1]{RenesseSturm:TransportIneq}, building work of Otto-Villani, \cite{OttoVillani:GenTalagrand}, and Cordero-Erausquin-McCann-Schmuckenschl\"ager in \cite{CMS:RiemInterpolationIneq}:
\begin{theo}\label{theo:riem-mfld-equiv}
For a compact Riemannian manifold $(M,g)$, regarding it as a measured geodesic space $\left(M, d, m\right)$, with $m=\widetilde{\vol_{M}}\in \cP(M)$, the normalized volume measure,\footnote{That is, we define $\widetilde{\vol_{M}}:= (\vol(M))^{-1} \vol_{M}$.} we have that $\Ent( \cdot | m)$ is weakly a.c.\ $K$-convex if and only if $M$ has the lower Ricci curvature bound,\footnote{By this, we mean as a bilinear form, or in other words $\Ric \geq K$ if and only if $\Ric(\xi,\xi') \geq Kg(\xi,\xi')$ for all $p\in M$ and $\xi,\xi' \in T_{p}M$.} $\Ric \geq K$ on $M$.\end{theo}
As such, we will often refer to a geodesic measure space $(X,d,\mu)$ for which $\Ent(\cdot|\mu)$ is weakly a.c.\ $K$-convex as a ``space with generalized $\Ric \geq K$''.

\subsection{Metric Structure of $(\cP_{0},d^{W})$} The principal reason that we work with measures on the unit interval is that the one dimensionality of the underlying space allows us use the inverse distribution function to embed $\cP_{0}:= \cP([0,1])$ as a convex subset of $L^{2}[0,1]$ whose Hilbert space structure proves crucial in proof of Theorem \ref{theo:QQbeta-no-ric-geqK} (we will remark on the difficulties present in extending our results to higher dimensions, where there is no such Hilbert space structure in Section \ref{sect:highdim}). 

\begin{prop}\label{prop:G0-P0-iso}
Letting $\cG_{0} \subset L^{2}([0,1])$ be the subset of the square integrable functions on the interval which are right continuous and nondecreasing as maps $g:[0,1]\to [0,1]$, and let $d^{L^{2}}$ be the metric induced on $\cG_{0}$ from the $L^{2}$ norm, then the map \begin{align*} \addlabel{eq:psi-defi} \Psi : (\cG_{0},d^{L^{2}}) &\to (\cP_{0},d^{W})\\
 g & \mapsto g_{*}\Leb \end{align*} is an isometry. The inverse $\Psi^{-1}$ is given by 
\[ \addlabel{eq:psi-inverse} \Psi^{-1} : \mu \mapsto g_{\mu}\] where $g_{\mu}$ is the inverse distribution function defined
\[ \addlabel{eq:inverse-dist-defi} g_{\mu}(s) : = \inf\{r \in[0,1]: \mu([0,r] )> s\}\] with the convention that $\inf \emptyset :=1$.
\end{prop} 

For a proof, see \cite[Theorem 2.18]{Villani:TopOptTspt03} or \cite[Proposition 4.1]{Chodosh:essay}. It is not hard to show that (see, for example \cite[Lemma 4.2]{Chodosh:essay}):
\begin{lemm}\label{lemm-g0-geo}
The space $\cG_{0}$ is a totally convex subset of $L^{2}([0,1])$ (that is, any geodesic between two elements in $\cG_{0}$ lies entirely in $\cG_{0}$). In fact, for $f,g\in \cG_{0}$, the unique geodesic between them is given by the linear combination
\[ \gamma(t) := (1-t) f + t g .\]
\end{lemm}
This fact proves highly beneficial to our analysis, and as we discuss in Section \ref{sect:highdim}, its failure to hold in higher dimensions is one of the obstacles in extending our results to $(\cP(M),d^{W},\PP^{\beta})$, for $M$ a general Riemannian manifold. As a result of the above proposition and lemma, instead of $(\cP_{0},d^{W})$, we can study the isometric space $(\cG_{0},d^{L^{2}})$, which is a convex subset of a Hilbert space, so has very simple geometry. 
\subsection{Entropic Measure}

We now discuss the measure $\QQ^{\beta}_{0} \in \cP(\cG_{0})$, originally constructed by von Renesse and Sturm in \cite[Proposition 3.4]{RenesseSturm:EntMeasDiff}. 
 \begin{defiprop}\label{defiprop-QQbeta}
 For $\beta >0$ there is a (unique) probability measure $\QQ^{\beta}_{0} \in \cP(\cG_{0})$ which we will call the \emph{entropic measure} (it could also be referred to as the law of the {Dirichlet process} or as the Gibbs measure) such that for each partition of $[0,1]$ \[ 0=t_{0}< t_{1}<\dots<t_{N}<t_{N+1} = 1\] and for all bounded measurable functions $u:[0,1]^{N} \to \RR$, we have that 
 \begin{multline}
 {\int_{\cG_{0}} u(g(t_{1}),\dots, g(t_{N}) )d \QQ^{\beta}_{0}}(g) \\
 = \frac{\Gamma(\beta)}{\prod_{i=0}^{N} \Gamma(\beta(t_{i+1}-t_{i}))} \int_{\Sigma_{N}} u(x_{1},\dots,x_{N}) \prod_{i=0}^{N} (x_{i+1} - x_{i})^{\beta(t_{i+1}-t_{i})-1} dx_{1}\cdots dx_{N}
 \end{multline}
 where we define \begin{equation*}\Sigma_N := \{ (x_1,\dots,x_N)\in [0,1]^N : 0 = x_0 < x_1 < \dots < x_N < x_{N+1} = 1\},\end{equation*} and $\Gamma(s) = \int_{0}^{\infty} t^{s-1}e^{-t}dt$ is the Gamma function. 
 \end{defiprop}
In their paper, von Renesse and Sturm prove the above existence result using the Kolmogorov extension theorem. We will denote the pushforward measure $\PP^{\beta}_{0} := \Psi_{*}\QQ^{\beta}_{0}$, where $\Psi$ is defined in Proposition \ref{prop:G0-P0-iso}. We are in reality interested in the metric measure triple $(\cP_{0},d^{W},\PP^{\beta}_{0})$, but because of its linear structure as a convex subset of $L^{2}$, we will find it far easier to work with the space $(\cG_{0},d^{L^{2}},\QQ^{\beta}_{0})$ (which is metric measure isomorphic to the space we are interested in, so provides an equivalent object for study).

\subsection{Heuristics} \label{subsect:heuristics}

In \cite{RenesseSturm:EntMeasDiff}, von Renesse and Sturm give a heuristic argument that $\PP^{\beta}_{0}$ is of the form 
\begin{equation}
d\PP^{\beta}_{0}(\mu) = \frac{1}{Z_{\beta}} e^{-\beta \Ent(\mu|m)} d\PP^{0}_{0}(\mu)
\end{equation}
where $Z_{\beta}$ is a normalizing constant, and $\Ent(\mu|m)$ is the entropy of $\mu$ with respect to the Lebesgue measure and $\PP^{0}_{0}$ is to be thought of as a ``uniform measure'' on $\cP_{0}$ (which does not actually exist). At first sight, one should expect such a measure to be displacement convex, because the $\Ent(\mu|m)$ is $0$-convex on $\cP_{0}$ by Theorem \ref{theo:riem-mfld-equiv}, and by \cite[Proposition 4.14]{Sturm:MetMeas1}, multiplying by the exponential of a convex function will not decrease generalized Ricci bounds, were they to exist for the hypothetical measure $\PP^{0}_{0}$. 

We also remark that the measure $\PP^{\beta}_{0}$ displays properties which are consistent with such lower bounds. Von Renesse and Sturm have constructed a symmetric Dirichlet form in \cite{RenesseSturm:EntMeasDiff}, $(\EE,D(\EE))$, given as the closure in $L^{2}(\cG_{0})$ of the quadratic form 
\[ 
\EE(F) : = \int_{\cG_{0}} |\DD F(g)|^{2}_{L^{2}([0,1])} d\QQ^{\beta}_{0}(g)
\]
with domain 
\[\left\{F(g) = \varphi \left(\bangle{f_{1},g}_{L^{2}([0,1])}, \dots, \bangle{f_{m},g}_{L^{2}([0,1])} \right) : m \geq 1, \varphi \in C^{1}_{b}(\RR^{m}), f_{k} \in L^{2} \right\}\]
and where $\DD F(g)$ is the $L^{2}$-Fr\'echet derivative of $F$ at $g$, which for $F$ in the domain described above is 
\[ 
\DD F(g) (x) = \sum_{i=1}^{m} \partial_{i}\varphi \left(\bangle{f_{1},g}_{L^{2}([0,1])}, \dots, \bangle{f_{m},g}_{L^{2}([0,1])} \right) f_{i}(x).
\]
The existence of such a Dirichlet form is interesting for various reasons (e.g.\ see \cite{FOT:DirichletForms}) but in our case, it is relevant because D\"oring and Stannat have shown that $\EE$ satisfies a Poincar\'e inequality
 \begin{theo}[\cite{Stannat:logsobWass} Theorem 1.2]\label{theo:poincineq}
The Dirichlet form constructed in \cite{RenesseSturm:EntMeasDiff}, $\EE$ satisfies a Poincar\'e inequality with constant less than $\frac 1\beta$, i.e.\ for all $F \in D(\EE)$
 \[ \Var_{\QQ^{\beta}_{0}}(F) \leq  \frac 1 \beta \EE(F). \]
 \end{theo}
 as well as a log-Sobolev inequality
 \begin{theo}[\cite{Stannat:logsobWass} Theorem 1.4] 
 There exists a constant $C$ (independent of $\beta$) such that for $F \in D(\EE)$
 \[ \int_{\cG_{0}} F(g)^{2} \log\left(\frac{F(g)^{2}}{ \| F \|^{2}_{L^{2}(\QQ^{\beta}_{0})}} \right) d\QQ^{\beta}_{0}(g)\leq \frac C \beta \EE(F).\]
 \end{theo}
Both of these theorems are properties that would hold true if $(\cG_{0},d^{L^{2}},\QQ^{\beta}_{0})$ had generalized Ricci curvature bounds (c.f. \cite[Corollary 6.12, Theorem 6.18]{LottVillani:RicciViaTspt}), (and in particular, Theorem \ref{theo:poincineq} would be the consequence of the space having generalized Ricci bounded below by $\beta$).\footnote{We remark that the log-Sobolev and Poincar\'e inequalities are certainly weaker conditions than generalized lower Ricci bounds. To see this, consider the possibility of such inequalities on a metric measure space $(X,d,m)$. It is clear that for a measurable function on $X$, $\rho(x)$ with $0<\alpha\leq \rho(x) \leq \beta <\infty$ and $\int \rho dm = 1$, the space $(X,d,\rho m)$ admits log-Sobolev and Poincar\'e inequalities if and only if $(X,d,m)$ does (with constants changing in a manner easily prescribed by $\alpha,\beta$. On the other hand, generalized Ricci bounds change (or are destroyed) in a much more sensitive manner depending on $\rho$, and should be thought of a ``higher order'' condition than log-Sobolev/Poincar\'e inequalities.  We thank the referee for drawing our attention to this point.  }

However, in spite of these heuristics, there are no such Ricci lower bounds, as we see in Theorem \ref{theo:QQbeta-no-ric-geqK}. We do remark that in \cite{Gigli:secondOrder}, Gigli has argued that there is no natural choice of volume form on $\cP(M)$, because this would be equivalent to there existing a Laplacian (by an integration by parts formula), which seems not to exist, because of the issues related to tracing a Hessian type object over an infinite dimensional space. In addition, he has written (in \cite[Remark 5.5]{Gigli:secondOrder}) that Sturm has communicated to him a measure theoretic argument that the measure $\PP^{\beta}_{0}$ could not be the volume form on $\cP_{0}$ (and the same for the higher dimensional analogue), and as such it seems we should stop searching for a volume measure on $\cP(M)$. Of course, none of this precludes an arbitrary reference measures in $\cP(\cP(M))$ giving rise to lower Ricci bounds on $\cP(M)$ (in particular, a point mass, say $\delta_{\text{Leb}} \in \cP(\cP([0,1]))$ certainly has such bounds, but is not a particularly interesting example).

\subsection{Acknowledgments}

This paper is a consolidated version of my essay written for Part III of the Mathematical Tripos at Cambridge University for the 2010-2011 academic year. I would like to thank Cl\'ement Mouhot for agreeing to set and mark the essay, assisting me in learning the material contained within, his extensive editing help, suggesting to me the problem of Ricci bounds on $(\cP(X),d^{W})$, as well as his continued support after the conclusion of the program. I would additionally like to thank Cedric Villani, as well as the anonymous referee for their helpful comments and suggestions. I am grateful to the Cambridge Gates Trust for their financial support during my year at Cambridge. Preparation of this paper was partially completed while supported by the National Science Foundation Graduate Research Fellowship under Grant No.\ DGE-1147470.

\section{Proof of Theorem \ref{theo:QQbeta-no-ric-geqK}} Suppose Theorem \ref{theo:QQbeta-no-ric-geqK} is false, so there is some $\beta> 0,K\in \RR$ such that the above space has generalized $\Ric \geq K$. We will show that this yields a contradiction, as follows. 

Certainly, without loss of generality, we may assume that $K \leq 0$. Let, for $s \in (0,1)$
\[ A_{s}: = \{g \in \cG_{0} : g(s) > 1/2 \}\]
\[ B_{s} := \{g \in\cG_{0} : g(s) > 0 \}.\]
It is clear that the following convex combination of these sets in $\cG_{0}$ (or equivalently $L^{2}$ by Lemma \ref{lemm-g0-geo}) is
\begin{equation} \label{eq:AB-convex-comb-C-def}(1-t) A_{s} + t B_{s} = \{ g\in \cG_{0} : g(s) > (1-t)/2\} := C_{s}(t). \end{equation} 
The significance of $C_{s}(t)$ is that any geodesic $\gamma:[0,1]\to \cG_{0}$ such that $\gamma(0) \in A_{s}$ and $\gamma(1) \in B_{s}$, has $\gamma(t) \in C_{s}(t)$. This will play a crucial role in our argument, allowing us to estimate the entropy of a Wasserstein geodesic in $\cP(X)$ between $\QQ^{\beta}_{0}$-uniform measures supported on $A_{s}$ and $B_{s}$.

Notice that 
\[ \QQ^{\beta}_{0}(C_{s}(t)) = \frac{\Gamma(\beta)}{\Gamma(\beta s) \Gamma(\beta (1-s))} \int_{(1-t)/2}^{1} x^{\beta s - 1}(1-x)^{\beta (1-s) -1} dx.\]
In particular, by Euler's beta integral (see \cite[Section 1.5]{Lebedev:SpecialFunctions}), we see that \[\QQ^{\beta}_{0}(B_{s}) = \QQ^{\beta}_{0}(C_{s}(1))  = 1,\] and for $s \in (0,1)$ and it is not hard to see that for all $t \in[0,1]$, $\QQ^{\beta}_{0} (C_{s}(t)) > 0$. As such, we, define 
\[ \mu({s}) := \frac{1}{\QQ^{\beta}_{0}(A_{s})} \chi_{A_{s}} \QQ^{\beta}_{0} \in \cP(\cG_{0}) \]
and by assumption, there is a geodesic $\mu(s)_{t}$ between $\mu(s)$ and $\QQ^{\beta}_{0}$ such that the entropy, $\Ent(\mu(s)_{t}|\QQ^{\beta}_{0})$ is $K$-convex, as in the definition of generalized $\Ric \geq K$.  Because \[\Ent(\mu(s) | \QQ^{\beta}) = \int_{\cG_{0}}  \frac{1}{\QQ^{\beta}_{0}(A_{s})} \chi_{A_{s}}  \log\left( \frac{1}{\QQ^{\beta}_{0}(A_{s})} \chi_{A_{s}} \right) d\QQ^{\beta}_{0} = - \log(\QQ^{\beta}_{0}(A_{s})) < \infty\]
(and clearly $\Ent(\QQ^{\beta}_{0}|\QQ^{\beta}_{0}) = 0$) so by assumptions of $K$-convexity, we must have that \[\Ent(\mu(s)_{t}|\QQ^{\beta}_{0}) <\infty,\] in particular implying that $\mu(s)_{t} \ll \QQ^{\beta}_{0}$. Thus, we can write $\mu(s)_{t} = \rho(s)_{t} \QQ^{\beta}_{0}$, and we have that by \eqref{eq:AB-convex-comb-C-def}, we see that $\mu(s)_{t}$ is concentrated in $C_{s}(t)$\footnote{This follows from the fact that optimal transport maps mass along geodesics of the underlying space. This is intuitively obvious, as if not, we could move along a geodesic between the endpoints, reducing the total distance traveled, and it follows rigorously from \cite[Proposition 2.10]{LottVillani:RicciViaTspt}. Thus, $\mu(s)_{t}$ is concentrated in $\cup_{\gamma} \gamma(t)$ where the union is over all geodesics $\gamma:[0,1] \to \cG_{0}$ with $\gamma(0) \in A_{s}$ and $\gamma(1) \in B$. Because $\cG_{0}$ is a totally geodesic subset of a Hilbert space, we have that this union is just $C_{s}(t)$. }, implying that
\begin{align*} \Ent(\mu(s)_{t}| \QQ^{\beta}_{0}) & = \int_{\cG_{0}} \rho(s)_{t}\log\rho(s)_{t} d\QQ^{\beta}_{0}\\
& = \int_{\cG_{0}} \log \rho(s)_{t} d\mu(s)_{t}\\
& = \int_{\cG_{0}} \log\left( \frac{d\mu(s)_{t}}{d\QQ^{\beta}_{0}|_{C_{s}(t)}} \frac{1}{\QQ^{\beta}_{0}(C_{s}(t)) }\right) d\mu(s)_{t}\\
& = \int_{\cG_{0}} \log\left( \frac{d\mu(s)_{t}}{d\QQ^{\beta}_{0}|_{C_{s}(t)}}\right) d\mu(s)_{t} -  \int_{\cG_{0}} \log\left({\QQ^{\beta}_{0}(C_{s}(t)) }\right) d\mu(s)_{t} \\
& = \Ent(\mu(s)_{t}|\QQ^{\beta}_{0}|_{C_{s}(t)}) - \log\left({\QQ^{\beta}_{0}(C_{s}(t)) }\right)\\
& \geq  - \log\left({\QQ^{\beta}_{0}(C_{s}(t)) }\right)
\end{align*}
where $\QQ^{\beta}_{0}|_{C_{s}(t)} = \frac{1}{\QQ^{\beta}_{0}(C_{s}(t))} \chi_{C_{s}(t)} \QQ^{\beta}_{0}$.
Combining this with the assumed $K$-convexity of the entropy functional along the path $\mu(s)_{t}$, we thus have that 
\begin{align*}
- \log \left(\QQ^{\beta}_{0}(C_{s}(t))\right) & \leq \Ent( \mu(s)_{t}|\QQ^{\beta}_{0})\\
&  \leq (1-t)\Ent(\mu(s)|\QQ^{\beta}_{0}) + t\underbrace{\Ent(\QQ^{\beta}_{0}|\QQ^{\beta}_{0})}_{=0} - \frac K 2 t(1-t) d(\mu(s),\QQ^{\beta}_{0})^{2}.
\end{align*}

This implies that, because $\Ent(\mu(s)|\QQ^{\beta}_{0}) = -\log\left( \QQ^{\beta}_{0}(A_{s})\right)$
\begin{equation}\label{eq:log-convex-gen-K} \log\left(\QQ^{\beta}_{0} (C_{s}(t)) \right) \geq (1-t) \log \left( \QQ^{\beta}_{0}(A_{s})\right) + \frac K 2 t(1-t) d^{W}(\mu(s),\QQ^{\beta}_{0})^{2} .\end{equation}
Because $\diam(\cG_{0}) = 1$, we must have that $d^{W}(\mu(s),\QQ^{\beta}) \in[0,1]$, so rearranging \eqref{eq:log-convex-gen-K} and using this observation, we have that (we've assumed that $K \leq 0$)
\begin{equation}
 \log \left( \frac{\QQ^{\beta}_{0}(C_{s}(t))}{(\QQ^{\beta}_{0}(A_{s}))^{1-t}} \right) \geq \frac K 2 t(1-t),
 \end{equation}
 implying that for all $s,t \in (0,1)$
 \begin{equation}\label{eq:final-cond-from-log-conc}
  \frac{\QQ^{\beta}_{0}(C_{s}(t))}{(\QQ^{\beta}_{0}(A_{s}))^{1-t}} \geq \exp\left({\frac K 2 t(1-t)}\right).
 \end{equation}
 We will show that for a fixed $t\in(0,1)$ 
 \begin{equation} \label{eq:lim-final-cond-from-log-conc}
  \lim_{s\to 0}   \frac{\QQ^{\beta}_{0}(C_{s}(t))}{(\QQ^{\beta}_{0}(A_{s}))^{1-t}} = 0,
 \end{equation}
  contradicting \eqref{eq:final-cond-from-log-conc}, because the right hand side is bounded away from zero for a fixed $t$ and $K\in \RR$. To see this, note that by definition of $\QQ^{\beta}_{0}$, in Definition-Proposition \ref{defiprop-QQbeta}, we have that
 \begin{equation}
  \frac{\QQ^{\beta}_{0}(C_{s}(t))}{(\QQ^{\beta}_{0}(A_{s}))^{1-t}}  = \left( \frac{\Gamma(\beta)}{\Gamma(\beta s)\Gamma(\beta (1-s))}\right)^{t} \frac{\int_{(1-t)/2}^{1} x^{\beta s - 1}(1-x)^{\beta(1-s) -1} dx}{\left( \int_{1/2}^{1} x^{\beta s -1} (1-x)^{\beta(1-s) -1} \right)^{1-t}}.
 \end{equation}
 It is not hard to see that because we have fixed $t\in(0,1)$, all of the terms have finite and nonzero limits (for this argument, it is important that we take $t < 1$ because the $x^{\beta s-1}$ term in the top integral becomes non-integrable at $x=0$ as $s\to 0$, but $t<1$ ensures that the limits of the integral do not include $0$), except for $\Gamma(\beta s)$ which approaches $\infty$, showing \eqref{eq:lim-final-cond-from-log-conc}, and thus showing that there cannot be any generalized Ricci lower bounds on $(\cG_{0},d^{L^{2}},\QQ^{\beta}_{0})$.

\section{Higher Dimensional Setting}\label{sect:highdim}

We remark that in \cite{Sturm:MultEnt}, Sturm has constructed a higher dimensional analogue of $\PP_{0}^{\beta}$, over a general compact Riemannian manifold, $M$. We briefly describe his construction and then explain why the method of proof of Theorem \ref{theo:QQbeta-no-ric-geqK} to $(\cP(M), d^{W},\PP^{\beta})$ does not seem to extend to this case. In order to discuss Sturm's construction, we first need the following:
\begin{defi}[$d^{2}/2$-concavity]
A function $\phi:M\to \RR$ is called \emph{${d^{2}}/{2}$-concave} if there exists a function $\psi: M\to \RR$ so that 
\[ \phi(y) =  \inf_{x\in M}\left[ \frac 12 d(x,y)^{2} - \psi(x) \right] \] for all $y \in M$. 
For a function $\psi: M\to \RR$, we define its \emph{${d^{2}}/{2}$-transform} $\phi^{\frac{d^{2}}{2}}$ by the same formula
\[ \phi^{\frac{d^{2}}{2}}(y) : = \inf_{x\in M} \left[\frac 12 d(x,y)^{2} - \phi(x)\right] .\]
\end{defi}

\begin{theo}[\text{\cite[Theorem 8]{McCann:PolarFactor}}] \label{theo:polar-fact-riem-mflds}  
For $\mu,\nu \in \cP(M)$ with $\mu \ll \vol_{M}$, there is a ${d^{2}}/{2}$-concave function $\phi:M\to \RR$ so that the map \[ F_{t} (x) := \exp_{x}(-t\nabla\phi)\] gives $\mu_{t} : = (F_{t})_{*}\mu : [0,1] \to \cP(M)$, which is the unique geodesic between $\mu$ and $\nu$. Furthermore, $(\Id,F_{1})_{*}\mu_{0}$ is an optimal transport plan between $\mu_{0}$ and $\mu_{1}$. For $t \in[0,1)$ $\mu_{t}\ll \vol_{M}$ and if, in addition, $\nu \ll \vol_{M}$, then we have that for all $t\in[0,1]$, $\mu_{t}\ll \vol_{M}$. 
\end{theo}

In fact, Sturm shows that $(\cP(M),d^{W})$ is homeomorphic to the space of maps of the form $g=\exp(-\nabla\phi)$ for $\phi$ a $d^{2}/2$-concave function $M\to \RR$ equipped with an $H^{1}$ Sobolev norm (the homeomorphism is given by $g\mapsto g_{*}\widetilde{\vol_{M}}$, i.e.\ pushing forward the normalized volume measure, $\widetilde{\vol_{M}}:= (\vol_{M}(M))^{-1} \vol_{M}$, and the inverse is given by finding the unique $d^{2}/2$-concave function associated to $\mu \in \cP(M)$ as given by Theorem \ref{theo:polar-fact-riem-mflds}). Sturm then notes that the $d^{2}/2$-transform map thus gives a involution $\cP(M) \to \cP(M)$, so for $\mu \in \cP(M)$ we map
\begin{equation} \mathfrak{C}: \mu\mapsto \exp(-\nabla\phi^{\frac{d^{2}}{2}})_{*}\widetilde{\vol_{M}}.  \end{equation}
Sturm then defines the measure $\QQ^{\beta}\in \cP(\cP(M))$ by requiring that for each measurable partition $M = \biguplus_{i=1}^{N} M_{i} $ and bounded Borel function $u:\RR^{N}\to \RR$, the following holds, where $m_{i} = \widetilde{\vol_{M}}(M_{i})$
\begin{multline}
\int_{\cP(M)} u(\nu(M_{1}),\dots,\nu(M_{N})) d\QQ^{\beta}(\nu) = \\
\frac{\Gamma(\beta)}{\prod_{i=1}^{N}\Gamma(\beta m_{i}) } \int_{[0,1]^{N}, \sum_{i=1}^{N} x_{i} = 1} u(x_{1},\dots,x_{N}) x_{1}^{\beta m_{1} -1  }\cdots x_{N}^{\beta   m_{N} -1  } dx_{1}\cdots dx_{N}.
\end{multline}
Sturm argues that this measure exists by the Kolomogorov extension theorem (and Euler's beta integral), and then defines the ``multidimensional entropic measure'' $\PP^{\beta} := \mathfrak{C}_{*}\QQ^{\beta}$.   

We conjecture that $(\cP(M), d^{W},\PP^{\beta})$ does not have an generalized lower Ricci bounds, but so far we have been unable to prove it. The principal difficulty seems to be that $\mathfrak{C}_{*}:\cP(\cP(M)) \to \cP(\cP(M))$ does not seem to map Wasserstein geodesics to geodesics, so that if we were to try to mimic the one dimensional case,  taking a geodesic from $\PP^{\beta}$ to a measure ``close'' to a singular measure, and then taking a sequence of such geodesics that end up closer and closer to the singular measure, to be able to compute the entropy we would have to push these geodesics forward by $\mathfrak{C}_{*}$, where they would no longer be geodesics. In particular, we remark that $\mathfrak{C}_{*}$ makes the measure $\PP^{\beta}$ into a very nonlocal object, in the sense that it seems that there is no result implying that if $M_{0} \subset M$ is a totally convex subset and we were to know that $(\cP(M),d^{W},\PP^{\beta})$ had generalized $\Ric\geq K$ then we could conclude that $(\cP(M_{0}), d^{W},\PP^{\beta}|_{M_{0}})$ (where $\PP^{\beta}|_{M_{0}}$ is the renormalization of $\PP^{\beta}$ restricted to the set of measures $\mu \in \cP(M)$ with $\mu(M_{0}) = 1$). This is because the $d^{2}/2$-transform seems to be quite sensitive to the geometry of the underlying space (for example the $d^{2}/2$-transform of $\delta_{x_{0}}$ is a singular measure which one can show is supported on the cut locus of $x_{0}$ and whose exact form seems to potentially have a complicated relationship with the geometry of $M$) which makes an argument similar to our proof of Theorem \ref{theo:QQbeta-no-ric-geqK} seem difficult in this higher dimensional setting.

\bibliography{bibliography}
\bibliographystyle{siam}
\end{document}